\documentclass[11pt]{amsart}
\usepackage[utf8x]{inputenc}
\usepackage{graphicx,color}

\newtheorem{thm0}{Theorem}[section]
\newtheorem{lemma}[thm0]{Lemma}
\newtheorem{prop}[thm0]{Proposition}
\newtheorem{rem}[thm0]{Remark}
\newtheorem{cor}[thm0]{Corollary}
\newtheorem{defin}[thm0]{Definition}
\newtheorem{exa}[thm0]{Example}

\def\rig#1{\smash{ \mathop{\longrightarrow}
    \limits^{#1}}}

\def\sear#1{\searrow
   \rlap{$\vcenter{\hbox{$\scriptstyle#1$}}$}}
\def\swar#1{\swarrow
   \rlap{$\vcenter{\hbox{$\scriptstyle#1$}}$}}
\def\dow#1{\Big\downarrow
   \rlap{$\vcenter{\hbox{$\scriptstyle#1$}}$}}

\def\O{{\mathcal O}}

\def\N{{\mathbb N}}

\def\C{{\mathbb C}}

\def\P#1{{\mathbb P}^#1}
\def\PP{{\mathbb P}}
\newcommand{\Vop}{{ V_0 \otimes \dots \otimes V_p}}
\newcommand{\Woq}{{ W_0 \otimes \dots \otimes W_q}}

\begin{document}

\title[Introduction to the hyperdeterminant]{Introduction to the hyperdeterminant and to the rank of multidimensional matrices}
\author{Giorgio Ottaviani}
\email{ottavian@math.unifi.it}
\address{Dipartimento di Matematica, Università di Firenze,  viale Morgagni 67/A, 50134 Firenze, ITALY }
\dedicatory{dedicated to David Eisenbud on the occasion of his $65$th birthday}
\maketitle 
\tableofcontents

\section{Introduction}

The classical theory of determinants was placed on a solid basis by Cayley in 1843. A few years later,
Cayley himself elaborated a generalization to the multidimensional setting \cite{Cay},
in two different ways. There are indeed several ways to generalize the notion of determinant to multidimensional matrices. Cayley's second
 attempt has a geometric flavour and was very fruitful.
This invariant constructed by Cayley is named today hyperdeterminant (after \cite{GKZ}) and reduces to the determinant in the case of square matrices,
which will be referred to as the classical case. The explicit computation
of the hyperdeterminant 
 presented from the very beginning exceptional difficulties . Even today explicit formulas
are known only in some cases, like the so called boundary format case and in a few others.
In general one has to invoke elimination theory.
Maybe for this reason the theory was forgotten for almost 150 years. 
Only in 1992, thanks to a fundamental paper by Gelfand, Kapranov and Zelevinsky,
the theory was placed in the modern language and many new results have been found.
The book \cite{GKZ}, of the same three authors, is the basic source on the topic. 
Also chapter 9 of \cite{W} is a recommended reading, a bit more advanced,
see also \cite{BW}.
Two sources about classical determinants are \cite{Muir} and \cite{Pas} (the second one has also a German translation). 
The extension of the determinant to the multidimensional setting contained in these two sources is based
on the formal extension of the formula computing the classical determinant summing over all the permutations
(like in Cayley first attempt), and they have different
 properties from the hyperdeterminant studied in \cite{GKZ}
(the paper \cite{Ghe} glimpses a link between the two approaches in the $2\times 2\times 2$ case).

In this survey we introduce the hyperdeterminants and some of its properties from scratch. Our aim is to provide elementary arguments,
when they are available. The main tools we use are the biduality theorem
and the language of vector bundles. We will use Geometric Invariant Theory only in Section \ref{git}.
Essentially no results are original,
but the presentation is more geometric than the standard one. In particular the basic computation of the dimension
of the dual to the Segre variety is performed by describing the contact locus in the Segre varieties.

I wish to thank an anonymous referee for careful reading and several useful suggestions.

\section{Multidimensional matrices and the local geometry of Segre varieties}

Let $V_i$ be complex vector spaces of dimension $k_i+1$ for $i=0,\ldots , p$.

We are interested in the tensor product $V_0\otimes\ldots\otimes V_p$,
where the group $GL(V_0)\times\ldots\times GL(V_p)$ acts in a natural way.

Once a basis is fixed in each $V_i$, the tensors can be represented
as multidimensional matrices of format
$(k_0+1)\times\ldots\times (k_p+1)$.

There are $p+1$ ways to cut a matrix of format
$(k_0+1)\times\ldots\times (k_p+1)$ into parallel slices,
generalizing the classical description of rows and columns for $p=1$.

\begin{figure}[h]
   \centering
\includegraphics[width=40mm]{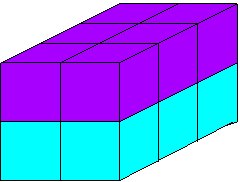}%
\qquad\qquad
\includegraphics[width=40mm]{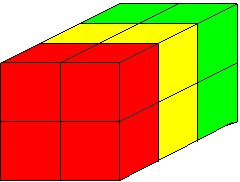}
    \caption{Two ways to cut a $3\times 2\times 2$ matrix into parallel slices}
\label{firstfigure}
    \end{figure}

The classical case $p=1$ is much easier than the case $p\ge 2$ mainly because there are only finitely many orbits for the action of $GL(V_0)\times GL(V_1)$. 

 Let \begin{equation}\label{defdr}
D_r=\{f\in V_0\otimes V_1|\textrm{rk\ }f\le r\}\end{equation}

We have that $D_r\setminus D_{r-1}$ are exactly the orbits of this action, and in particular the maximal rank matrices form the dense orbit.

Note that $D_1$ is isomorphic to the Segre variety $\PP(V_0)\times\PP(V_1)$ (they were introduced in \cite{Se})
and that it coincides with the set of decomposable tensors,
which have the form $v_0\otimes v_1$ for $v_i\in V_i$.

The first remark is 

\begin{lemma}\label{secrank1} The rank of $f$ coincides with the minimum number $r$ of summands
in a decomposition $f=\sum_{i=1}^rt_i$ with $t_i\in D_1$.
\end{lemma}
\proof Acting with the group $GL(V_0)\times GL(V_1)$
$f$ takes the form $f=\sum_{i=1}^rv_0^i\otimes v_1^i$, where $\{v_0^i\}$ is a basis of $V_0$
and $\{v_1^i\}$ is a basis of $V_1$,
corresponding to the matrix

$$\left[\begin{array}{cc}
I_r&0\\
0&0\\
\end{array}\right]
$$

In this form the statement is obvious.\qed
\vskip 0.5cm

The $k$-th secant variety $\sigma_k(X)$ of a projective irreducible
variety $X$ is the Zariski closure of the union of the projective span
$<x_1,\ldots x_k>$ where $x_i\in X$.
We have a chain of inclusions

$$X=\sigma_1(X)\subset\sigma_2(X)\subset\ldots$$
With this definition, Lemma \ref{secrank1} reads
\begin{cor}\label{sigmakd1}
$$\sigma_k(D_1)=D_k$$
\end{cor}
Let's state also, for future reference, the celebrated ``Terracini Lemma'' (see e.g. \cite{Zak}), whose proof is
straightforward by a local computation.
\begin{thm0}[Terracini Lemma]\label{terracini}
Let $X$ be a projective irreducible variety and let $z\in<x_1,\ldots, x_k>$ be a general point in $\sigma_k(X)$.
Then $$T_z\sigma_k(X)=<T_{x_1}X,\ldots, T_{x_k}X>$$
\end{thm0}
The tangent spaces $T_{x_i}X$ appearing in the Terracini lemma are the projective 
tangent spaces. Sometimes, we will denote by the same symbol the affine tangent
spaces, this abuse of notation should not create any serious confusion.
\vskip 0.5cm
We illustrate a few properties of the Segre variety $\P{{}}(V_0)\times\ldots\times\P{{}}(V_p)$.

It is, in a natural way, a projective variety according to the Segre embedding
$$\begin{array}{ccc}\P{{}}(V_0)\times\ldots\times\P{{}}(V_p)&\rig{}&\P{{}}(V_0\otimes\ldots\otimes V_p)\\
(v_0,\ldots v_p)&\mapsto &v_0\otimes\ldots\otimes v_p\end{array}$$

In this embedding, the Segre variety
coincides with the projectivization of the set of {\it decomposable} tensors.
The proof of the following proposition is straightforward (by induction on $p$) and we omit it.
\begin{prop} \label{rankcontraction}Every $\phi\in V_0\otimes\ldots\otimes V_p$ induces, for any $i=0,\ldots,p$ the contraction map
$$C_i(\phi)\colon V_1^{\vee}\otimes\ldots \widehat{V_i^{\vee}}\ldots\otimes V_p^{\vee}\rig{}V_i$$
where the $i$-th factor is dropped from the source space.
The tensor $\phi$ is decomposable if and only if rk$(C_{i}(\phi))\le 1$ for every $i=0,\ldots, p$.
\end{prop}

\begin{figure}[h]
   \centering
\includegraphics[width=30mm]{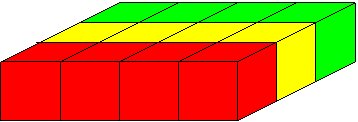}\quad
\includegraphics[width=30mm]{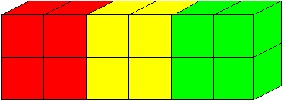}\quad
\includegraphics[width=30mm]{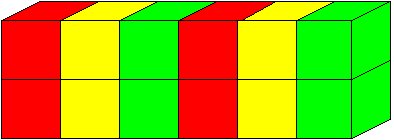}
    \caption{The three flattenings of the matrix in Figure \ref{firstfigure}. If the $2$-minors of two of them vanish, 
then the matrix corresponds to a decomposable tensor (a point in the Segre variety)}
\label{figflattening}
    \end{figure}

The previous proposition gives equations of the Segre variety as $2\times 2$ minors of the contraction maps
$C_i(\phi)$. These maps are called flattenings, because they are represented by bidimensional matrices
obtained like in Figure \ref{figflattening}.

\begin{rem} In Prop. \ref{rankcontraction} it is enough that the rank conditions are satisfied for all $i=0,\ldots,p$ except one .
\end{rem}

A feature of the Segre variety is that it contains a lot of linear subspaces.

For any point  $x=v_0\otimes\ldots\otimes v_p$,
the linear space  $v_0\otimes\ldots V_i\ldots\otimes v_p$ 
passes through $x$ for $i=0,\ldots,p$; it can be identified
with the fiber of the projection
$$\pi_i\colon\P{{k_0}}\times\ldots\times\P{{k_p}}\rig{}\P{{k_0}}\times\ldots\widehat{\P{{k_i}}}\ldots\times\P{{k_p}}$$
We will denote the projectivization
of the  linear subspace $v_0\otimes\ldots V_i\ldots\otimes v_p$ 
 as $\P{{k_i}}_x$.

These linear spaces have important properties described by the following proposition.

\begin{prop}\label{tangentsegre}
Let $x\in X=\P{{k_0}}\times\ldots\times\P{{k_p}}$.

(i) The tangent space at $x$ is the span of the $p+1$ linear spaces $\P{{k_i}}_x$ ,
that is $T_pX$ is the projectivization of $\oplus_i v_0\otimes\ldots V_i\ldots\otimes v_p$

(ii) The tangent space at $x$ meets $X$ in the union of the $p+1$ linear spaces $\P{{k_i}}_x$.

(iii) Any linear space in $X$ passing through $x$ is contained in one of the $p+1$ linear spaces $\P{{k_i}}_x$.
\end{prop}

\proof The tangent vector to a path $v_0(t)\otimes\ldots\otimes v_p(t)$ for $t=0$ is
$\sum_{i=0}^pv_0(0)\otimes\ldots v_i'(0)\ldots \otimes v_p(0)$. Since $v_i'(0)$
may be chosen as an arbitrary vector, the statement (i) is clear.

(ii) Fix a basis $\{e_j^0,\ldots , e_j^{k_j}\}$ of $V_j$ for $j=0,\ldots ,p$
and let $\{e_{j,0},\ldots , e_{j,k_j}\}$ be the dual basis. We may assume that $x$ corresponds to 
 $e_0^0\otimes\ldots\otimes e_p^0$. Consider a decomposable tensor $\phi$ in the tangent space at $x$, so
$\phi=v_0\otimes\ldots\otimes e_p^0+\ldots+e_0^0\otimes\ldots\otimes v_p$ for some $v_i$.
We want to prove that  $v_i$ and $e_i^0$ are linearly independent for at most one index $i$.
Otherwise we may assume $\dim(v_0,e_0^0)=2$, $\dim(v_1,e_1^0)=2$. 
Consider the contraction
 $$\begin{array}{ccccc}C_0(\phi)(e_{1,0}\otimes e_{2,0}\otimes\ldots\otimes e_{p,0})&=&v_0&+&(\ldots)e_0^0\\
C_0(\phi)(e_{1,1}\otimes e_{2,0}\otimes\ldots\otimes e_{p,0})&=&&&\left(e_{1,1}(v_1)+\sum_{j=2}^pe_{j,0}(v_j)\right)e_0^0
\end{array}$$
Since we may assume also $e_{1,1}(v_1)\neq 0$, by replacing $e_{1,1}$ with a scalar multiple we have also
$\left(e_{1,1}(v_1)+\sum_{j=2}^pe_{j,0}(v_j)\right)\neq 0$. This implies that $\textrm{rank\ }C_0(\phi)\ge 2$
which is a contradiction. For an alternative approach generalizable to any homogeneous space see \cite{LM}.

(iii) A linear space in $X$ passing through $x$ is contained in the tangent space at $x$,
hence the statement follows from (ii).\qed

\section{The biduality Theorem and the contact loci in the Segre varieties}

The projective space ${\PP}(V)$ consists of linear subspaces of dimension one of $V$.
The dual space  ${\PP}(V^\vee)$ consists of linear subspaces of codimension one (hyperplanes) of $V$.
Hence the points in ${\PP}(V^\vee)$ are exactly the hyperplanes of  ${\PP}(V)$. 

Let's recall the definition of dual variety. Let $X\subset {\PP}(V)$ be a projective irreducible variety  . A hyperplane $H$ is called
{\it tangent} to  $X$ if $H$ contains the tangent space to  $X$ at some nonsingular point $x\in X$.

The {\it dual variety} $X^{\vee}\subset {\PP}(V^\vee)$  is defined as the Zariski closure of the set of all the
tangent hyperplanes. Part of the biduality theorem below says that ${X^{\vee}}^{\vee}=X$,
but more is true.
Consider the incidence variety $V$ given by the closure of the set
$$\{(x,H)\in X\times {\PP}(V^\vee)|x\textrm{\ is a smooth point and\ }T_xX\subset H\}$$
 $V$ is identified in a natural way with the projective bundle ${\PP}(N(-1)^{\vee})$
where $N$ is the normal bundle to $X$ (see Remark \ref{normal}).

\begin{thm0}(Biduality Theorem)
Let $X\subset {\PP}(V)$ be an irreducible projective variety. We have
\begin{equation}
{X^{\vee}}^{\vee}=X  
\end{equation}
Moreover if $x$ is a smooth point of $X$ and $H$ is a smooth point of $X^{\vee}$,
then $H$ is tangent to $X$ at $x$ if and only if $x$, regarded as a hyperplane in ${\PP}(V^\vee)$,
is tangent to $X^{\vee}$ at $H$.
In other words the diagram
\begin{equation}\label{dualdiagram}\begin{array}{ccccc}
&&V\\
&\swar{p_1}&&\sear{p_2}\\
X&&&&X^{\vee}\\
\end{array}\end{equation}
is symmetric.
\end{thm0}

For a proof, in the setting of symplectic geometry, we refer to \cite{GKZ}, Theorem 1.1 .

Note, as a consequence of the biduality theorem, that the fibers of both the projections of $V$ over smooth points are linear spaces.
This is trivial for the left projection, but it is not trivial for the right one.
Let's record this fact

\begin{cor}\label{contactlinear}
Let $X$ be smooth and let $H$ be a general tangent hyperplane (corresponding
to a smooth point of $X^{\vee}$).
Then $\{x\in X| T_xX\subseteq H\}$ is a linear subspace (this is called the contact locus of $H$ in $X$).
\end{cor}

As a first application we compute the dimension of the dual to a Segre variety.

\begin{thm0}\label{dimdual}[Contact loci in Segre varieties]
Let $X={\PP}^{k_0}\times \ldots \times {\PP}^{k_p}$.

(i) If $k_0\ge \sum_{i=1}^pk_i$ then a general hyperplane tangent at $x$ is tangent along a linear space
of dimension $k_0-\sum_{i=1}^pk_i$ contained in the fiber $\P{{k_0}}_x$ .
In this case the codimension of $X^{\vee}$ 
is $1+k_0-\sum_{i=1}^pk_i$.

(ii) If $k_0\le \sum_{i=1}^pk_i$ then a general hyperplane tangent at $x$ is tangent only at $x$.
In this case $X^{\vee}$  is a hypersurface.

(iii) The dual variety  $X^{\vee}$
is a hypersurface if and only if the following holds 
 $$\max k_i=k_0\le\sum_{i=1}^pk_i$$
\end{thm0}

\proof We remind that, by Proposition \ref{tangentsegre} (i),
a hyperplane $H$ is tangent at $x$ if and only if it contains the $p+1$ fibers 
through $x$. By Corollary \ref{contactlinear} a general hyperplane is tangent along a linear variety. By 
Prop. \ref{tangentsegre} (iii) a linear variety in $X$ it is contained in one of the fibers.
Let $H$ be a general hyperplane tangent at $x$.
We inspect the fibers through $y$ when $y\in \P{{k_0}}_x$. The locus where $H$ contains
the fiber $\P{{k_i}}_y$  is a linear space in $\P{{k_0}}_x$ of codimension $k_i$,
indeed the fibers can be globally parametrized by $y$ plus other $k_i$ independent points.
This description
 proves (i), because the variety $V$ in (\ref{dualdiagram}) has the same dimension of a hypersurface in $\PP(V^{\vee}$ and 
we just computed the general fibers of $p_2$. Also (ii) follows by the same argument because the conditions
are more than the dimension of the space.
(iii) is a consequence of (i) and (ii).  \qed
\begin{figure}[h]
   \centering
\includegraphics[width=40mm]{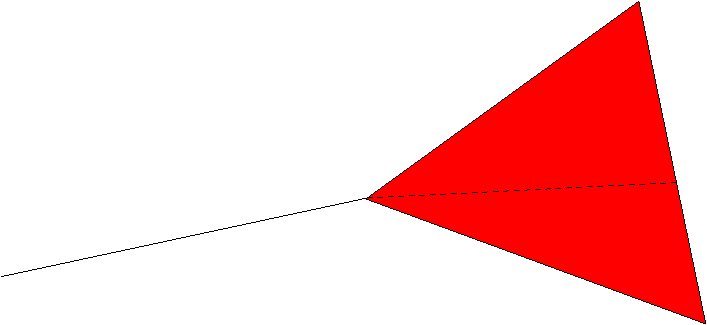}
    \caption{The tangent space at a point $x\in X=\PP^1\times\PP^2$ cuts $X$ into two linear spaces meeting at $x$,
the general hyperplane tangent at $x$ is tangent along a line (dotted in the figure)}
    \end{figure}
\begin{defin}\label{defboundary}
A format $(k_0+1)\times\ldots\times(k_p+1)$ with $k_0=\max_j{k_j}$ is called a
boundary format if $k_0=\sum_{i=1}^pk_i$.
In other words, the boundary format 
corresponds to the equality in (iii) of Theorem \ref{dimdual}
\end{defin}
\begin{rem}\label{normal} According to [L], Theorem \ref{dimdual} says that for a Segre variety with normal bundle $N$, the twist
$N(-1)$ is ample if and only if the inequality $\max k_i=k_0\le\sum_{i=1}^pk_i$ holds.
\end{rem}
Note that for $p=1$ the dual variety to $D_1={\PP}^{k_0}\times {\PP}^{k_1}$
is a hypersurface if and only if $k_0=k_1$ (square case). This is better understood
by the following result.

\begin{thm0}\label{dualp1}Let $k_0\ge k_1$. In the projective spaces of $(k_0+1)\times (k_1+1)$ matrices
the dual variety to the variety $D_r$ (defined in formula (\ref{defdr})) is $D_{k_1+1-r}$.

When $k_0=k_1$ (square case) the determinant hypersurface is the dual of $D_1$.
\end{thm0}

In order to prove the Theorem \ref{dualp1} we need the following proposition
\begin{prop}\label{inclusionsecant}
Let $X$ be a irreducible projective variety. For any $k$ 
$$\left(\sigma_{k+1}(X)\right)^{\vee}\subset \left(\sigma_{k}(X)\right)^{\vee}$$
\end{prop}
\proof The proposition is almost a tautology after the Terracini Lemma.
The dual to the $(k+1)$-th secant variety $\left(\sigma_{k+1}(X)\right)^{\vee}$ is defined as the closure of the set of
 hyperplanes $H$ containing $T_{z}\sigma_{k+1}(X)$ for $z$ being a smooth point in $\sigma_{k+1}(X)$,
so $z\in <x_1,\ldots ,x_{k+1}>$ for general $x_i\in X$. 
By the Terracini Lemma (Prop. \ref{terracini})
$H$ contains $T_{x_1},\ldots T_{x_{k+1}}$, hence  
 $H$ contains
$T_{z'}\sigma_{k}(X)$ for the general $z'\in<x_1,\ldots ,x_{k}>$
(removing the last point).
\qed

{\it Proof of Theorem \ref{dualp1}}.
Due to Proposition \ref{inclusionsecant} and Corollary \ref{sigmakd1} we have the chain of inclusions
$$D_1^{\vee}\supset D_2^{\vee}\supset\ldots\supset D_{k_1}^{\vee}$$
By the biduality Theorem any inclusion must be strict.
Since the $D_i^{\vee}$ are $GL(V_0)\times GL(V_1)$-invariant, 
and the finitely many orbit closures are given by $D_i$, the only possible solution
is that the above chain coincides with
$$D_{k_1}\supset\ldots\supset D_1$$.\qed

\begin{exa} When $X\subset \P n$ is the rational normal curve,
$\sigma_{k}(X)$ consists of polynomials which are sums of $k$ powers,
while  $\sigma_{k}(X)^{\vee}$
consists of polynomials having $k$ double roots. We get that $\sigma_{k}(X)^{\vee}=Chow_{2^k,1^{n-2k}}(\P 1)$
according to the notations of Section \ref{symmetric}.
\end{exa}

\begin{rem} A common misunderstanding after Theorem \ref{dualp1} is that $X\subset Y$ implies
the converse inclusion $X^{\vee}\supset Y^{\vee}$. This is in general false. The simplest counterexample
is to take $X$ to be a point of a smooth (plane) conic $Y$. Here $X^{\vee}$ is a line and $Y^{\vee}$ is again a smooth conic.
\end{rem}

\begin{rem}
The proof of Theorem \ref{dualp1} is short, avoiding local computations, but rather indirect.

 We point out the elegant proof of Theorem \ref{dualp1} given by Eisenbud in
Prop. 1.7 of \cite{E88}, which gives more information. Eisenbud considers $V_0\otimes V_1$ as the space of linear maps
$Hom(V_0^{\vee},V_1)$ and its dual $Hom(V_1,V_0^{\vee})$. These spaces are dual under the pairing
$<f,g>:=tr(fg)$ for $f\in Hom(V_0^{\vee},V_1)$ and $g\in Hom(V_1,V_0^{\vee})$. Eisenbud proves that
if $f\in D_r\setminus D_{r-1}$ then the tangent hyperplanes at $f$ to $D_r$ are exactly
the $g$ such that $fg=0$, $gf=0$. These conditions force
the rank of $g$ to be $\le k_1+1-r$. Conversely any $g$ of rank $\le k_1+r-1$
satisfies these two conditions for some $f$ of rank $r$, proving Theorem \ref{dualp1} .
\end{rem}

The above proposition is
 important because it gives a geometric interpretation of the determinant, as the dual
of the Segre variety. This is the notion that better generalizes to multidimensional matrices.

\begin{defin} \label{defhyper}
Let
 $$\max k_i=k_0\le\sum_{i=1}^pk_i$$
The equation of the dual variety to  ${\PP}^{k_0}\times \ldots \times {\PP}^{k_p}$
 is the hyperdeterminant.
\end{defin}

A point deserves a clarification. Since the dual variety lives in the dual space,
we have defined the hyperdeterminant in the dual space to the space of matrices,
and not in the original space of matrices.
Although there is no canonical isomorphism between the space of matrices and its dual space
this apparent ambiguity can be solved by the invariance.

Indeed we may construct infinitely many isomorphisms between the space

$V_0\otimes\ldots\otimes V_p$ and its dual  $V_0^{\vee}\otimes\ldots\otimes V_p^{\vee}$,
by fixing a basis constructed from the bases of the single spaces $V_i$.
Any function on the space of matrices which is invariant with respect to the action
of  $SL(V_0)\times\ldots\times SL(V_p)$
induces, by using any isomorphisms $V_i\simeq V_i^{\vee}$, a function on the dual space
and, due to the invariance, it does not depend on the chosen isomorphism.

\section{Degenerate matrices and the hyperdeterminant}
From now on we will refer to a multidimensional matrix $A$ as simply a matrix,
and write $Det(A)$ for its hyperdeterminant (when it exists).
\begin{defin}\label{defindegenere}A matrix $A$ is called degenerate if
there exists a nonzero $ (x^0\otimes x^1\otimes\ldots\otimes x^p)\in V_0\otimes \ldots\otimes V_p$  such that 
\begin{equation}\label{formuladefdegenere}
A(x^0,x^1,\ldots ,V_i,\ldots, x^p)=0\qquad\forall i=0,\ldots , p 
\end{equation}

\end{defin}

The ``kernel'' $K(A)$ is by definition the variety of  nonzero $ (x^0\otimes x^1\otimes\ldots\otimes x^p)\in V_0\otimes \ldots\otimes V_p$ such that 
(\ref{formuladefdegenere}) is satisfied.

So, by Proposition \ref{tangentsegre} (i),  a matrix $A$ is degenerate if and only if $A$ corresponds to an hyperplane tangent in $K(A)$. This gives an algebraic reformulation of the definition of hyperdeterminant.

We get
\begin{prop}\label{degequaldual}
(i) The (projectivization of the) variety of degenerate matrices of format $(k_0+1)\times\ldots\times (k_p+1)$
is the dual variety of the Segre variety ${\PP}^{k_0}\times\ldots\times {\PP}^{k_p}$.

(ii) Let $\max k_i=k_0\le\sum_{i=1}^pk_i$. A matrix $A$ is degenerate if and only if 
$Det (A)=0$.
\end{prop}
\proof
Part (i) is a reformulation of  Proposition \ref{tangentsegre} (i) .
Part (ii) is a reformulation of Theorem \ref{dimdual} (iii) and Definition \ref{defhyper}. \qed

Degenerate matrices $A$ such that $K(A)$ consists of a single point are exactly the smooth points
of the hypersurface.

\begin{figure}[h]
   \centering
\includegraphics[width=40mm]{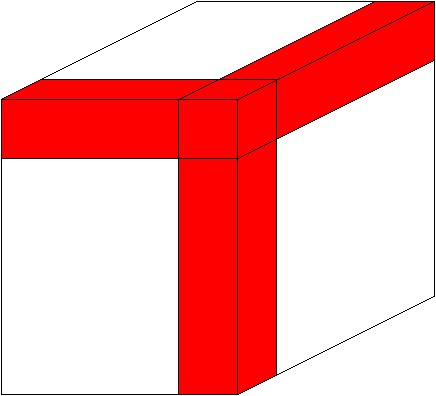}
    \caption{A matrix $A$ is degenerate if and only if it becomes zero on the colored part after a linear change of coordinates.
The kernel $K(A)$ corresponds to the colored vertex.}
\label{coloredvertex}
    \end{figure}
\vskip 0.5cm

The degree of the hyperdeterminant can be computed by means of a generating function. Let  $N(k_0,\ldots ,k_p)$ 
be the degree of the hyperdeterminant of format  $(k_0+1)\times\ldots\times(k_p+1)$ and let $k_0=\max_{j}k_j$ .
Set $N(k_0,\ldots ,k_p)=0$ if $k_0>\sum_{i=1}^pk_i$ (in this case we may set  $Det =1$).
\begin{thm0}[\cite{GKZ} Thm. XIV 2.4]
\[\sum_{k_0,\ldots, k_p\ge 0}N(k_0,\ldots ,k_p)z_0^{k_0}\cdots z_p^{k_p}
=\frac{1}{\left(1-\sum_{i=2}^{p+1}(i-1)x_i(z_0,\ldots, z_p)\right)^2}\]
where $x_i$ is the $i$-th elementary symmetric function.
\end{thm0}

For example for $p=2$,
$$\sum_{k_0,k_1, k_2\ge 0}N(k_0,k_1,k_2)z_0^{k_0}z_1^{k_1}z_2^{k_2}
=\frac{1}{\left(1-(z_0z_1+z_0z_2+z_1z_2)-2z_0z_1z_2\right)^2}$$

We list a few useful degree of hyperdeterminants,
corresponding to formats $(a,b,c)$ with $a+b+c\le 12$

$$\begin{array}{c|c|c}
\textrm{format}&\textrm{degree}&\textrm{boundary format}\\
\hline\\
(2,2,2)&4\\
(2,2,3)&6&*\\
(2,3,3)&12\\
(2,3,4)&12&*\\
(2,4,4)&24\\
(2,4,5)&20\\
(3,3,3)&36\\
(3,3,4)&48\\
(3,3,5)&30&*\\
(3,4,4)&108\\
(3,4,5)&120\\
(4,4,4)&272\\
(2,b,b)&2b(b-1)&\\
(2,b,b+1)&b(b+1)&*\\
(a,b,a+b-1)&\frac{(a+b-1)!}{(a-1)!(b-1)!}&*\\
\end{array}$$

Note that the degree of format $(2,b,b+1)$ is smaller than the degree of its subformat $(2,b,b)$ (for $b\ge 4$).
Therefore a Laplace expansion cannot exist, at least not in a naive way.

In the boundary format case the degree simplifies to $\frac{(k_0+1)!}{k_1!\ldots k_p|}$,
as we will see in Section \ref{multilinear}.

Let's see what happens to the hyperdeterminant after swapping two parallel slices.

\begin{thm0}
Let $N$ be the degree of hyperdeterminant of format  $(k_0+1)\times\ldots\times(k_p+1)$ .
(i) $\frac{N}{k_i+1}$ is an integer.
(ii) After swapping two parallel slices of format $(k_0+1)\times\ldots\widehat{(k_i+1)}\ldots\times(k_p+1)$
the hyperdeterminant changes its sign if $\frac{N}{k_i+1}$ is odd and remains invariant if
$\frac{N}{k_i+1}$ is even.

(iii) A matrix with two proportional parallel slices has hyperdeterminant equal to zero.
\end{thm0}
\proof It is clear from its definition that
the hyperdeterminant is a relative invariant for the group
$G=GL(k_0+1)\times\ldots\times GL(k_p+1)$. Moreover the hyperdeterminant is homogeneous on each slice,
and by the action, the degree has to be the same on parallel slices. Since there are $(k_i+1)$ parallel slices,
the hyperdeterminant is homogeneous of degree $\frac{ N}{(k_i+1)}$ with respect to each slice,
which proves (i).
Hence for $g\in GL(k_i+1)$ we get 
$Det(A*g)=Det(A)\cdot(\det(g))^{N/(k_i+1)}$. If $g$ is a permutation matrix,
it acts on the parallel slices by permuting them. In particular if $g$ swaps two slices it satisfies  $\det g=-1$,
hence (ii) follows . (iii) follows because a convenient $g$ acting on slices
makes a whole slice equal to zero.\qed

\begin{figure}[h]
    \centering
\includegraphics[width=40mm]{iniziocubo2_0.jpg}%
\qquad\qquad$\Longrightarrow$\qquad\qquad
\includegraphics[width=40mm]{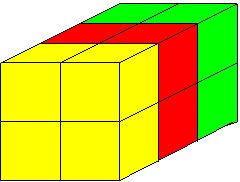}
    \caption{Swapping two ``vertical'' slices in the format  $3\times 2\times 2$ leaves the hyperdeterminant invariant.}
\end{figure}

\begin{figure}[h]
    \centering
\includegraphics[width=40mm]{iniziocubo3_2.jpg}%
\qquad\qquad$\Longrightarrow$\qquad\qquad
\includegraphics[width=40mm]{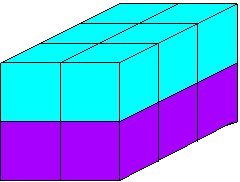}
    \caption{Swapping two ``horizontal'' slices in the format  $3\times 2\times 2$ 
the hyperdeterminant changes sign.}
\end{figure}
 
\section{Schl\"afli technique of computation}

For simplicity we develop Schl\"afli technique in this Section
only for three-dimensional matrices, although a similar argument works in any dimension (see \cite{GKZ}).

Let $A$ be a matrix of format $a\times b\times b$. It can be seen as a $b\times b$ matrix 
with entries which are linear forms over $V_0$ and we denote it by
$\tilde A(x)$ with $x\in V_0\simeq\C^{a}$ .

\begin{figure}[h]
    \centering
\includegraphics[width=25mm]{iniziocubo2_0.jpg}
\qquad$\Longrightarrow$\qquad\qquad
$z_0$\includegraphics[width=15mm]{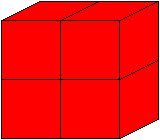}$\quad + z_1$\includegraphics[width=15mm]{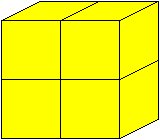}
$\quad + z_2$\includegraphics[width=15mm]{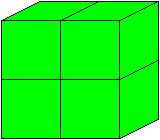}
    \caption{A tridimensional matrix gives a bidimensional matrix with linear forms as coefficients.}
\end{figure}

Then we can compute  $\det\tilde A(x)$
and we get a homogeneous form over $V_0$ that we can identify with a hypersurface in $\PP{V_0}$.

Schl\"afli main remark is the following.

\begin{thm0}[Schl\"afli]\label{schlafli}

Let $A$ be degenerate and let $v_0\otimes v_1\otimes v_2\in K(A)$.

Then the hypersurface $\det\tilde A(x)$ is singular at $v_0$.
\end{thm0}
\proof After a linear change of coordinates we may assume that
$v_0\otimes v_1\otimes v_2$ corresponds to the corner of the matrix $A$,
meaning that we choose a basis of $V_i$ where $v_i$ is the first element.
Let $\tilde A=\sum_{i=0}^{k_1}x_i A_i$
where $A_i$ are the $b\times b$ slices.
After a look at Figure \ref{coloredvertex}, we get that
the assumption that $v_0\otimes v_1\otimes v_2\in K(A)$ means that the first slice has the form 

$$ A_0=\left[\begin{array}{ccccc}
0&0&\ldots,&0\\
0&*&\ldots&*\\
0&\vdots&&\vdots\\
0&*&\ldots&*\\
\end{array}\right]$$
and the successive slices have the form

$$ A_i=\left[\begin{array}{ccccc}
0&*&\ldots,&*\\
*&*&\ldots&*\\
*&\vdots&&\vdots\\
*&*&\ldots&*\\
\end{array}\right]$$

Note that $\det A_0=0$
and $v_0$ has coordinates $(1,0,\ldots 0)$.
In an affine coordinate system centered at this point we may assume $x_0=1$.

In conclusion we get 

$$\tilde A=\left[\begin{array}{ccccc}
0&\sum_{i}l_{1,i}x_i&\ldots,&\sum_{i}l_{k_1,i}x_i\\
\sum_{i}m_{1,i}x_i&*&\ldots&*\\
\vdots&\vdots&&\vdots\\
\sum_{i}m_{k_1,i}x_i&*&\ldots&*\\
\end{array}\right]$$

where $l_{j,i}$ and $m_{j,i}$ are scalars.

Now, by expanding the determinant on the first row and then on the first column, we get
that $\det\tilde A$ has no linear terms in $x_i$
and then the origin is a singular point, as we wanted.\qed

The conclusion of Theorem \ref{schlafli} is that for a degenerate matrix $A$, the discriminant of the polynomial $\det\tilde A(x)$ has to vanish. In other words the hyperdeterminant of $A$ divides the discriminant of $\det\tilde A(x)$.
The following proposition characterizes exactly the cases where the converse holds and it was used by Schl\"afli to
compute the hyperdeterminant in these cases.

\begin{prop}\label{4bb}

(i) For matrices $A$ of format $2\times b\times b$ and $3\times b\times b$  the discriminant of $\det\tilde A(x)$ coincides with the hyperdeterminant of $A$.

(ii) For matrices $A$ of format $4\times b\times b$ the discriminant of $\det\tilde A(x)$ is equal to the product of the hyperdeterminant and an extra factor
which is a square.

(iii) For matrices $A$ of format $a\times b\times b$ with $a\ge 5$, the discriminant of $\det\tilde A(x)$ vanishes identically. 
\end{prop}

{\it Partial proof} (see \cite{GKZ} for a complete proof) (iii) follows because the locus of matrices of rank $\le b-2$ has codimension $4$ in the space
of $b\times b$ matrices, then every determinant hypersurface in a space of projective dimension at least
four contains points where the rank drops at least by $2$, and these are singular points of the hypersurface.
In particular for $a\ge 5$, $\det\tilde A(x)$ is singular.

Let's prove (i) in the format case $2\times b\times b$ .
From Theorem \ref{schlafli} we get that the hyperdeterminant of $A$ divides the discriminant of $\det\tilde A(x)$.
The degree of the discriminant of a polynomial of degree $b$ is $2(b-1)$, hence the degree of 
the discriminant of $\det\tilde A(x)$ is $2b(b-1)$ with respect to the coefficients of $A$.
Once we know that the degree of the hyperdeterminant of format $2\times b\times b$
is $2b(b-1)$ the proof is completed. This was stated in the table in the previous section
but it was not proved. We provide an alternative argument
for $\det\tilde A(x)$ singular implies $A$ is degenerate. As a consequence this provides also a proof of the degree formula.
After a linear change of coordinates we may assume as in the proof of Theorem \ref{schlafli} that 
$$A_0=\left[\begin{array}{cc}
0&0\\
0&B_0\\
\end{array}\right]$$
 and 

$$ A_1=\left[\begin{array}{ccccc}
c_1&*\\
*&*\\
\end{array}\right]$$

Note that $\det A_0=0$
and $v_0$ has coordinates $(1,0,\ldots 0)$.
In an affine coordinate system centered at this point we may assume $x_0=1$.

Then $\det\tilde A=\det\left[\begin{array}{ccccc}
c_1x_1&x_1\cdot(*)\\
x_1\cdot(*)&B_0+x_1\cdot(*)\\
\end{array}\right]=c_1x_1\det(B_0)+(\textrm{higher terms in\ } x_i)$.
If 
$\det(B_0)\neq 0$ we get $c_1=0$, hence $A$ is degenerate. 
If $\det(B_0)=0$ we can rearrange $A_0$ in such a way that
the first two rows and two columns vanish. In this case it is easy to arrange a zero in the
upper left $2\times 2$ block of $A_1$,
so that again $A$ is degenerate.
This pattern can be extended to the $3\times b\times b$ case. \qed

\begin{exa} Regarding (ii) of Proposition \ref{4bb}, we describe the format $4\times 3\times 3$.
From the table in the previous section, the hyperdeterminant of format $4\times 3\times 3$ has degree $48$. The discriminant
of cubic surfaces has degree $32$, and $32\cdot 3=96$. So there is an extra factor
of degree $48$. It is the square of the invariant coming from the invariant $I_8$,
which is the invariant of minimal degree for cubic surfaces \cite{Dol}.

This gives an interpretation of the invariant $I_8$. Namely, a smooth cubic surface
has $I_8$ vanishing if and only if it has a determinantal representation
as a $4\times 3\times 3$ matrix which is degenerate.
\end{exa}

\begin{exa}\label{first322}
Let's see the $3\times 2\times 2$ example,
which can be computed explicitly in at least three different ways.
The first way is an application of Schl\"afli technique and it is described here. The second way is by looking at
multilinear systems and it is described in Theorem \ref{cayleypro}.
The third way is described in Example \ref{third}.

A matrix $A$ of format $3\times 2\times 2$ defines the following $2\times 2$ 
matrix
\[\left[\begin{array}{cc}
a_{000}x_0+a_{100}x_1+a_{200}x_2&a_{001}x_0+a_{101}x_1+a_{201}x_2\\
a_{010}x_0+a_{110}x_1+a_{210}x_2&a_{011}x_0+a_{111}x_1+a_{211}x_2\\
\end{array}\right]\]
The determinant of this matrix defines the following projective conic in the variables
$x_0, x_1, x_2$
\begin{equation}
\label{sch322}
x_0^2
\det\left[\begin{array}{cc}a_{000}&a_{001}\\
a_{010}&a_{011}\\
\end{array}\right]+
x_0x_1\left(\det\left[\begin{array}{cc}a_{000}&a_{101}\\
a_{010}&a_{111}\\
\end{array}\right]+
\det\left[\begin{array}{cc}a_{100}&a_{001}\\
a_{110}&a_{011}\\
\end{array}\right]\right)+\ldots =\end{equation}
\begin{equation}
=(x_0, x_1, x_2)\cdot C\cdot (x_0, x_1, x_2)^t
\end{equation}
where $C$ is a $3\times 3$ symmetric matrix.

By Proposition \ref{4bb}
the hyperdeterminant of $A$ 
is equal to the determinant of $C$, hence the hyperdeterminant vanishes
if and only if the previous conic is singular.

\end{exa}

\begin{exa}\label{333} An interesting example is given by matrices of format $3\times 3\times 3$.
The hyperdeterminant of format $3\times 3\times 3$ has degree $36$. The discriminant
of cubic curves has degree $12$, and $12\cdot 3=36$.

The three determinants obtained by the three possible directions give three elliptic cubic curves. 
When one of the three is smooth, then all three are smooth 
(this happens if and only if $Det\neq 0$) and all three are isomorphic
(see \cite{TC} and also \cite{Ng} Proposition 1). 
So we get three different determinantal representations of the same cubic curve. The Theorem 1 of \cite{Ng} says that the moduli space of 
$3\times 3\times 3$ matrices under the action of $GL(3)\times GL(3)\times GL(3)$
is isomorphic to the moduli space of triples $(C, L_1, L_2)$ where
$C$ is an elliptic curve and $L_1$, $L_2$ are two (non isomorphic) line bundles of degree three induced by
by the pullback of $\O(1)$ under the other two determinantal representations.
\end{exa}
\begin{exa}
The $2\times 2\times 2$ case gives with an analogous computation the celebrated
Cayley formula for $A$ of format $2\times 2\times 2$
from the discriminant of the polynomial

\[\det\left[\begin{array}{cc}
a_{000}x_0+a_{010}x_1& a_{001}x_0+a_{011}x_1\\
a_{100}x_0+a_{110}x_1& a_{101}x_0+a_{111}x_1\\
\end{array}\right]=0\]

which is
\[Det (A)=\left(\left|\begin{array}{cc}a_{000}&a_{011}\\
a_{100}&a_{111}\\
\end{array}\right|+
\left|\begin{array}{cc}a_{010}&a_{001}\\
a_{110}&a_{101}\\
\end{array}\right|\right)
^2-4 \left|\begin{array}{cc}a_{000}&a_{001}\\
a_{100}&a_{101}\\
\end{array}\right|
\cdot
\left|\begin{array}{cc}a_{010}&a_{011}\\
a_{110}&a_{111}\\
\end{array}\right|\]

The previous formula expands with exactly 12 summands which have a nice symmetry
and are the following

\cite{Cay}
{\footnotesize $$
Det(A)=\left(a_{000}^2a_{111}^2+a_{001}^2a_{110}^2+a_{010}^2a_{101}^2+a_{011}^2a_{100}^2\right)+$$
$$-2\left(a_{000}a_{001}a_{110}a_{111}+a_{000}a_{010}a_{101}a_{111}+
a_{000}a_{011}a_{100}a_{111}+\right.$$
$$\left.+a_{001}a_{010}a_{101}a_{110}+
a_{001}a_{011}a_{110}a_{100}+
a_{010}a_{011}a_{101}a_{100}\right)+$$
$$+4\left(a_{000}a_{011}a_{101}a_{110}+
a_{001}a_{010}a_{100}a_{111}\right)
$$
}

The first four summands correspond to the four line diagonals of the cube,
the following six summands correspond to the six plane diagonals while the last two summands correspond
to the vertices of the even and odd tetrahedra inscribed in the cube.

The $2\times 2\times 2$ hyperdeterminant is homogeneous of degree two in each slice
and remains invariant under swapping two slices. 
The $2\times 2\times 2$ hyperdeterminant appeared several times recently in the physics literature, see for example \cite{Du}.
\end{exa}
\begin{rem} According to \cite{HSYY}, the hyperdeterminant of format $2\times 2\times 2\times 2$
is a polynomial of degree $24$ containing $2,894,276$ terms. A more concise expression has been found in \cite{SZ},
in the setting of algebraic statistics,
by expressing the hyperdeterminant in terms of cumulants. The resulting expression has $13,819$ terms.
\end{rem}
\section{Multilinear systems, matrices of boundary format}
\label{multilinear}
A square matrix $A$ is degenerate if and only if the linear system $A\cdot x=0$
has a nonzero solution. This section explores how this notion generalizes to the multidimensional setting,
by replacing the linear system with a multilinear system (we borrowed some extracts from \cite{OV},
and I thank J. Vallès for his permission). 
The answer is that the hyperdeterminant captures the condition of existence of nontrivial solutions
only in the boundary format case.

Let $k_0=\max_j\{k_j\}$. A matrix $A$ of format $(k_0+1)\times\ldots\times(k_p+1)$
defines the linear map $C_0(A)\in Hom(V_1^{\vee}\otimes\ldots V_p^{\vee},V_0)$
(see Prop. \ref{rankcontraction})
which in turn defines a multilinear system. A nontrivial solution
of this system is given by nonzero $x^i\in V_i^{\vee}$ such that
$$C_0(A)(x^1\otimes\ldots\otimes x^p)=0$$

This is equivalent to the case $i=0$ of the definition of degenerate
(see (\ref{formuladefdegenere}) in Def. \ref{defindegenere}), namely to
\begin{equation}\label{condi0eq}
\exists\textrm{\ nonzero\ }x^1\otimes\ldots\otimes x^p\in V_1\otimes\ldots\otimes V_p
\textrm{\ such that\ }A(V_0,x^1,\ldots , x^p)=0\end{equation}

We say that $A$ satisfying (\ref{condi0eq}) is $0$-degenerate.

In the language of \cite{E88}, the $0$-degenerate matrices correspond exactly to the matrices
which are not $1$-generic. The condition to be not $0$-degenerate can be expressed indeed as a Chow condition
imposing that the linear subspace $\ker A$ meets the Segre variety.
For computations of Chow conditions see \cite{ESW}.

\begin{thm0}\label{condi0} [\cite{GKZ}, theor. XIV.3.1]
Let $k_0\ge\sum_{i=1}^pk_i$ (in particular in the boundary format). $A$ is $0$-degenerate if and only if it is degenerate.
\end{thm0}

\proof If $A$ is degenerate it is obviously $0$-degenerate. Let assume conversely
that $A$ is $0$-degenerate. By assumption there is a nonzero $(x^1\otimes\ldots\otimes x^p)\in V_1\otimes \ldots\otimes V_p$
such that (\ref{condi0eq}) holds.

Consider in the unknown $x^0$ the linear system (which is (\ref{formuladefdegenere}) for $i=1$)
$$A(x^0,V_1,x^2,\ldots , x^p)=0$$
It consists of $k_1$ independent equations with respect to (\ref{condi0eq}),
because one of the equations is already contained in (\ref{condi0eq}).
For the same reason the linear system given by  (\ref{formuladefdegenere})  for general $i$
consists of $k_i$ independent equations with respect to (\ref{condi0eq}).
All together we have $\sum_{i=1}^pk_i$ linear equations
in the $k_0+1$ unknowns which are the coordinates of $x^0$,
since the unknowns are more than the number of equations, we get a nonzero solution as we wanted. \qed

\begin{thm0}[\cite{GKZ} theor. XIV 1.3] \label{distriangolare}(Triangular inequality). 

(i) If $k_0\ge\sum_{i=1}^pk_i$ the variety of $0$-degenerate matrices has codimension $1+k_0-\sum_{i=1}^pk_i$ in $M(k_0+1,\ldots, k_p+1)$.

(ii) If $k_0<\sum_{i=1}^pk_i$ all matrices are $0$-degenerate.

(iii) The variety of $0$-degenerate matrices has codimension $1$ in $M(k_0+1,\ldots, k_p+1)$
exactly when the equality  $k_0=\sum_{i=1}^pk_i$ holds, that is in the boundary format case.
\end{thm0}

\proof If $k_0\ge\sum_{i=1}^pk_i$ the codimension is the same
as that of the variety of degenerate matrices by Theorem \ref{condi0}.
This codimension has been computed in Theorem \ref{dimdual} and Prop. \ref{degequaldual}
and it is $1$ only when the equality holds. 

If $k_0<\sum_{i=1}^pk_i$ all matrices are $0$-degenerate. Indeed
the kernel of $A\in Hom(V_1^{\vee}\otimes\ldots V_p^{\vee},V_0)$ meets the Segre variety by dimensional reasons.\qed

\begin{figure}[h]
   \centering
\includegraphics[width=40mm]{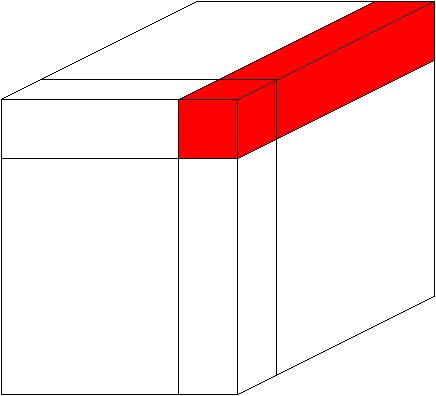}
    \caption{$A$ is $0$-degenerate if it becomes zero on the colored part after a linear change of 
coordinates}
    \end{figure}

We get the second promised expression for the  $3\times 2\times 2$ case.

\begin{thm0}\label{cayleypro} [Cayley] Let $A$ be a matrix of format $3\times 2\times 2$ and let $A_{00}$, $A_{01}$, 
  $A_{10}$, $A_{11}$ be the  $3\times 3$ submatrices obtained from
\[\left[\begin{array}{cccc}
a_{000}&a_{001}&a_{010}&a_{011}\\
a_{100}&a_{101}&a_{110}&a_{111}\\
a_{200}&a_{201}&a_{210}&a_{211}\\
\end{array}\right]\]
by removing respectively the first column $(00)$, the second $(01)$, the third $(10)$
and the fourth $(11)$.
The multilinear system $A(x\otimes y)=0$ given by  $\sum a_{ikj}x_jy_k=0$ has nontrivial solutions
if and only if 
$\det A_{01}\det A_{10}-\det A_{00}\det A_{11}=0$
which coincides with $Det (A)$.
\end{thm0}

{\it Proof\ } We may assume that the $3\times 4$ matrix has rank $3$,
otherwise we get a system with only two independent equations which has always a nontrivial solution.
A solution $(x_0y_0, x_0y_1, x_1y_0, x_1y_1)$ has to be proportional to $\left( \det A_{00}, -\det A_{01}, \det A_{10}, -\det A_{11} \right)$
by the Cramer rule. Now the condition is just the equation of the Segre quadric $\P1\times\P1\subset\P3$.\qed
\vskip 0.5cm

Let $A\in V_0\otimes\ldots\otimes V_p$ be of boundary format and let $m_j=\sum_{i=1}^{j-1}k_i$ with
the convention $m_1=0$.

We remark that the definition of $m_i$ depends on the order we have
chosen among the $k_j$'s (see Remark \ref{perm}).

With the above notations the vector spaces $V_0^\vee \otimes
S^{m_1}V_1\otimes\ldots\otimes S^{m_p}V_p$ and $
S^{m_1+1}V_1\otimes\ldots\otimes S^{m_p+1}V_p$ have the same dimension
$N=\frac{(k_0+1)!}{k_1!\ldots k_r!}.$

\begin{thm0}{\bf (and definition of $\partial_A$).}
\label{main}
Let $k_0=\sum_{i=1}^pk_i$. Then the hypersurface of $0$-degenerate matrices
has degree $N=\frac{(k_0+1)!}{k_1!\ldots k_r!}$ and its
equation is given by the determinant of the natural morphism
\[\partial_A: V_0^\vee \otimes S^{m_1}V_1\otimes\ldots\otimes
S^{m_p}V_p \rig{}
S^{m_1+1}V_1\otimes\ldots\otimes S^{m_p+1}V_p\]
\end{thm0}

\begin{proof}
  If $A$ is $0$-degenerate then we get $A(v_1\otimes\ldots\otimes v_p)=0$
  for some $v_i\in V_i^\vee$, $v_i\neq 0$ for $i=1,\ldots ,p$. Then
  $({\partial}_A)^{\vee}
\left(v_1^{\otimes m_1+1} \otimes\ldots\otimes v_p^{\otimes
      m_p+1} \right)=0$.
  
  Conversely if $A$ is non $0$-degenerate we get a surjective natural map
  of vector bundles over $X=\PP(V_2)\times\ldots\times\PP(V_p)$
\begin{equation*}
  \label{map}
 V_0^\vee \otimes \O_X\rig{\phi_A}V_1\otimes\O_X(1,\ldots ,1).  
\end{equation*}

Indeed, by our definition, $\phi_A$ is surjective if and only if $A$ is non $0$-degenerate.

We construct a vector bundle $S$ over
$\PP(V_2)\times\ldots\times\PP(V_p)$ whose dual $S^\vee$ is the kernel of
$\phi_A$ so that we have the exact sequence
\begin{equation}
  \label{ker}
0\rig{}S^\vee\rig{}V_0^\vee\otimes \O\rig{}V_1\otimes\O(1,\ldots ,1)\rig{}0. 
\end{equation}

After tensoring by $\O(m_2,\ldots ,m_p)$ and taking cohomology we get
\begin{equation*}
  \label{eq:3}
H^0(S^\vee(m_2,m_3,\ldots ,m_p))\rig{}V_0^\vee\otimes S^{m_1}V_1\otimes\ldots\otimes   S^{m_p}V_p\rig{{\partial}_A}
S^{m_1+1}V_1\otimes\ldots\otimes S^{m_p+1}V_p
\end{equation*}

and we need to prove
\begin{equation}
  \label{van}
H^0(S^\vee(m_2,m_3,\ldots ,m_p))=0.
\end{equation}

Let
$d=\dim\left(\PP(V_2)\times\ldots\times\PP(V_p)\right)=\sum_{i=2}^pk_i=m_{p+1}-k_1$.

Since $\det (S^\vee)=\O(-k_1-1,\ldots ,-k_1-1)$ and $rank~S^\vee=d$ ,
 it follows, by using the natural identification $S^{\vee}\simeq\wedge^{d-1}S\otimes \det (S^{\vee})$, that
$$S^\vee(m_2,m_3,\ldots ,m_p)\simeq\wedge^{d-1}S\otimes \det (S^\vee)(m_2,m_3,\ldots ,m_p)$$
hence
\begin{equation}
  \label{isos}
S^\vee(m_2,m_3,\ldots ,m_p)\simeq
\wedge^{d-1}S(-1,-k_1-1+m_3,\ldots,
-k_1-1+m_p) 
\end{equation}

Hence, by taking  the  $(d-1)$-st wedge power of the dual of
the sequence (\ref{ker}), and using  K\"unneth's formula 
to calculate the cohomology as in \cite{GKZ1}, the result follows.
\end{proof}

\begin{cor}
\label{defdet}
  Let $k_0=\sum_{i=1}^pk_i$. The hyperdeterminant of $A\in V_0\otimes\ldots\otimes V_p$ is the usual
  determinant of ${\partial}_A$, that is
  \begin{equation}
    \label{eq:partial}
Det(A):=det {\partial}_A    
  \end{equation}
where ${\partial}_A=H^0(\phi_A)$ and $\phi_A:V_0^\vee \otimes
\O_X\rig{\phi_A}V_1\otimes\O_X(1,\ldots ,1)$ is the sheaf morphism
associated to $A$. In particular 
\[deg Det =\frac{(k_0+1)!}{k_1!\ldots k_r!}\] 
\end{cor} 
This is  theorem 3.3 of chapter 14 of \cite{GKZ}.

\begin{rem}
\label{perm}
  \rm Any permutation of the $p$ numbers $k_1, \ldots, k_p$ gives
different $m_i$'s and hence a
  different map ${\partial}_A$. As noticed by Gelfand, Kapranov and 
Zelevinsky, in all cases the determinant of ${\partial}_A$ is the
  same by  Theorem \ref{main}. 
\end{rem}  

\begin{exa}\label{third} {\bf The $3\times 2\times 2$ case (third computation).}
In this case the morphism
$V_0^\vee\otimes V_1\to S^2V_1\otimes V_2$ is represented by a 
$6\times 6$ matrix, which, with the obvious notations, is the following

\begin{equation}
\label{sym}
\left[\begin{array}{cccccc}
a_{000}&a_{100}&a_{200}&0&0&0\\
a_{001}&a_{101}&a_{201}&0&0&0\\
a_{010}&a_{110}&a_{210}&a_{000}&a_{100}&a_{200}\\
a_{011}&a_{111}&a_{211}&a_{001}&a_{101}&a_{201}\\
0&0&0&a_{010}&a_{110}&a_{210}\\
0&0&0&a_{011}&a_{111}&a_{211}\\
\end{array}\right]
\end{equation}

The hyperdeterminant is the determinant of this matrix.
Note that this determinant is symmetric with respect to the second and the third index,
but this is not apparent from the above matrix.

\end{exa}

\begin{exa}
In the case $4\times 3\times 2$ the hyperdeterminant can be obtained
as the usual determinant of one of the following two maps
\[V_0^\vee\otimes V_1\to S^2V_1\otimes V_2\]
\[V_0^\vee\otimes S^2V_2\to V_1\otimes S^3V_2\]
\end{exa}
\begin{rem} The fact that  the degree of the hypersurface of $0$-degenerate
 matrices is
$N=\frac{(k_0+1)!}{k_1!\ldots k_p!}$ could be obtained in an alternative way.
We know 
 that A is $0$-degenerate iff
the corresponding $\ker A$ meets the Segre variety.

Hence the condition is given by a polynomial $P(x_1,\ldots ,x_m)$ in the 
variables
 $x_i\in\PP\left(\wedge^{k_0+1}(V_1\otimes\ldots\otimes V_p)\right)$ of degree equal to 
the degree of the Segre variety which is 
$\frac{k_0!}{k_1!\ldots k_p!}$. Since $x_i$ have degree $k_0+1$
in terms of the coefficients of $A$, the result follows.

\end{rem}
\vskip 0.5cm

Let $A=(a_{i_0, \dots, i_p})$ a matrix of boundary format $(k_0+1)\times \dots
\times (k_p +1)$ and $B=(b_{j_0, \dots, j_q})$ of boundary format
$(l_0+1)\times \dots \times (l_q +1)$, if $k_p=l_0$ the convolution  (or product)
 $A \ast B$ (see
\cite{GKZ}) of $A$ and $B$ is defined as the
$(p+q)$-dimensional matrix $C$ of format $(k_0+1)\times \dots \times
(k_{p-1} +1)\times(l_1+1)\times \dots \times (l_q +1)$ with entries
\begin{equation*}
  c_{i_{0},\dots,i_{p-1},j_1,\dots,j_q}=\sum_{h=0}^{k_p}a_{i_0, 
\dots,i_{p-1},h}b_{h,j_1,\dots, j_q}.
\end{equation*}
Note that $C$ has again boundary format. The following analogue of the
{\it Cauchy-Binet formula} holds.

\begin{thm0}
  \label{binet}
  If $A\in \Vop$ and $B\in \Woq$ are  boundary format
  matrices with $dim V_i=k_i+1$, $dim W_j=l_j+1$ and $W_0^\vee \simeq
  V_p$ then $A\ast B$ satisfies
\begin{equation}
  \label{eq:binet}
Det(A\ast
B)=(DetA)^{\binom{l_0}{l_1,\dots,l_q}}(DetB)^{\binom{k_0+1}{k_1,\dots,k_{p-1},k_p+1}}
\end{equation}
\end{thm0}
\begin{proof} 
\cite{DO}  
\end{proof}

\begin{cor}
$A$ and $B$ are nondegenerate if and only if $A*B$ is nondegenerate.
\end{cor}

\begin{exa}
  From Corolllary  \ref{defdet} the degree
  of the hyperdeterminant of a boundary format $(k_0+1)\times \dots
  \times (k_p+1)$ matrix $A$  is given by
   
  $$
  N_A=\frac{(k_0+1)!} 
  {k_1!\dots k_p!}
$$
Thus, (\ref{eq:binet}) can be rewritten as 
\begin{equation*}
  Det(A\ast B)={\lbrack (DetA)^{N_B} (DetB)^{N_A} \rbrack}^{\frac{1}{l_0+1}}
\end{equation*}
\end{exa}

\section{Link with Geometric Invariant Theory in the boundary format case}
\label{git}
\begin{figure}[h]
    \centering
\includegraphics[width=40mm]{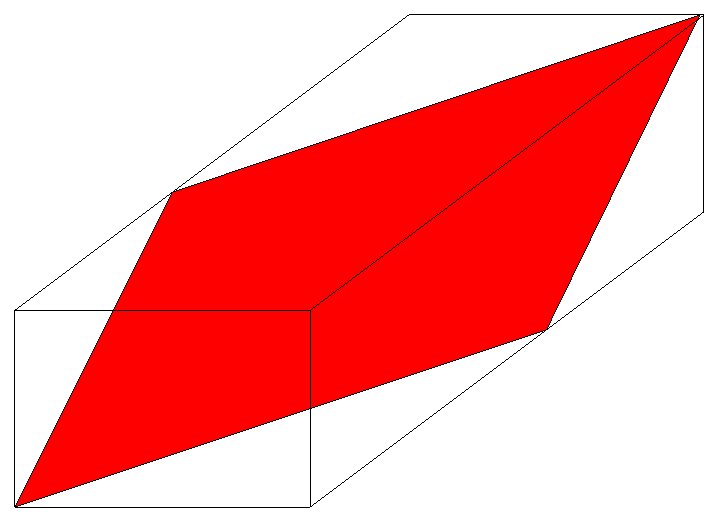}%
    \caption{The diagonal of a boundary format matrix. A triangular matrix fits one of the two half-spaces cut by the diagonal.}
\label{dia2}
\end{figure}

In the boundary format case it is well defined a unique ``diagonal''
given by elements $a_{i_0\ldots i_p}$ satisfying $i_0=\sum_{j=1}^pi_j$
(see Figure \ref{dia2}).
We will see in this section how these matrices behave under the action
of $SL(V_0)\times\ldots\times SL(V_p)$ in the setting of
the Geometric Invariant Theory.

\begin{defin} \label{x.1} A $p+1$-dimensional matrix of boundary
format $A\in V_0\otimes\ldots\otimes V_p$
 is  called  triangulable if   there exist  bases  in $V_j$ such that $a_{i_0,\ldots ,i_p}=0$
for $i_0>\sum_{t=1}^p i_t$

\end{defin}

\begin{defin} \label{x.2} A $p+1$-dimensional matrix of boundary
format $A\in V_0\otimes\ldots\otimes V_p$   is called
diagonalizable if   there exist  bases in $V_j$ such that $a_{i_0,\ldots ,i_p}=0$
for $i_0\neq\sum_{t=1}^p i_t$

\end{defin}
 
\begin{defin}\label{x.3}  A $p+1$-dimensional matrix of boundary format
 $A\in V_0\otimes\ldots\otimes V_p$ is  an identity if
one of the following equivalent
 conditions holds 
 
\noindent{i)} there exist  bases in $V_j$ such that   
$$a_{i_0,\ldots ,i_p}=\left\{\begin{array}{c}0\quad\hbox{for}\quad
i_0\neq\sum_{t=1}^pi_t\\
1\quad\hbox{for}\quad i_0=\sum_{t=1}^pi_t\\ \end{array}\right.$$

\noindent{ii)} there exist  a vector space $U$ of dimension $2$  and
isomorphisms $V_j\simeq S^{k_j}U$ such that
$A$ belongs to the unique one dimensional $SL(U)$-invariant subspace of
$S^{k_0}U\otimes S^{k_1}U\otimes\ldots\otimes S^{k_p}U$
\end{defin}

The equivalence between i) and ii) follows easily from the following remark:
 the matrix  $A$ satisfies the condition ii) if and only if
it corresponds to the natural multiplication map
$S^{k_1}U\otimes\ldots\otimes S^{k_p}U\to S^{k_0}U$ (after a suitable
isomorphism
$U\simeq U^\vee$ has been fixed).

The definitions of triangulable,
diagonalizable and identity
apply    to elements of  $\PP (V_0\otimes\ldots\otimes V_p)$ as well. In
particular
all identity matrices fill a distinguished orbit in
$\PP (V_0\otimes\ldots\otimes V_p)$. 

The function $Det$ is  $SL(V_0)\times\ldots\times SL(V_p)$-invariant, in
particular if $Det~A\neq 0$ then $A$ is semistable for the action of
$SL(V_0)\times\ldots\times SL(V_p)$.
We denote by $Stab~(A)\subset SL(V_0)\times\ldots\times SL(V_p)$
the  stabilizer subgroup of $A$ and by $Stab~(A)^0$ its connected component
containing
the identity. The main results  are the following.

\begin{thm0}\label{t1} ([AO]) Let $A\in \PP (V_0\otimes\ldots\otimes V_p)$
 of boundary format such that $Det~A\neq 0$. Then
$$A\hbox{\ is triangulable}\iff A\hbox{\ is not stable for the action of\ }
SL(V_0)\times\ldots\times SL(V_p)$$
\end{thm0}

\begin{thm0}\label{t2} ([AO]) Let $A\in \PP (V_0\otimes\ldots\otimes V_p)$
be of boundary format such that $Det~A\neq 0$.  Then
$$A\hbox{\ is diagonalizable}\iff Stab(A)\hbox{\ contains a subgroup isomorphic 
to\ }
\C^*$$
\end{thm0}

The proof of the above two theorems relies on the Hilbert-Mumford criterion.
The proof of the following theorem needs more geometry.

\begin{thm0} \label{t3} ([AO] for $p=2$, [D] for $p\ge 3$) Let $A\in \PP (V_0\otimes V_1
\otimes\ldots\otimes V_p)$
of boundary format such that $Det~A\neq 0$.  Then there exists a $2$-dimensional
vector space $U$ such that $SL(U)$ acts over $V_i\simeq S^{k_i}U$ and
according to this action on $V_0\otimes\ldots \otimes V_p$ we have
$Stab~(A)^0\subset SL(U)$. Moreover the following cases are possible
$$Stab~(A)^0\simeq\left\{\begin{array}{cl}0&\hbox{(trivial subgroup)}\\ \C\\ \C^*&\\
SL(2)&\hbox{\ (this case occurs if and only if\ }A\hbox{\ is an
identity)}\end{array}\right.$$
\end{thm0}

\begin{rem} When $A$ is an  identity then $Stab~(A)\simeq SL(2)$.
\end{rem}

\begin{exa}From the expression we have seen in the $3\times 2\times 2$ case one can compute that  the
hyperdeterminant of a diagonal matrix is

\[Det (A)=a_{000}^2a_{110}a_{101}a_{211}^2\]

In general the hyperdeterminant of a diagonal matrix is given by a monomial involving all the coefficients
on the diagonal
with certain exponents, see \cite{WZ}.
\end{exa}

\section{The symmetric case }\label{symmetric}

We analyze  the classical {\it pencil of quadrics} from the point of view of hyperdeterminants.

Let $A$ be a $2\times n\times n$ matrix
with two symmetric $n\times n$ slices $A_0$, $A_1$.
\begin{prop}\label{degensim} If $A$ is degenerate, then its kernel contains an element of the form
$z\otimes x\otimes x\in\C^2\otimes\C^n\otimes\C^n$.
\end{prop}
\proof By assumption there are nonzero $z=(z_0,z_1)^t\in\C^2$ and $x, y\in\C^n$ such that
$x^tA_iy=0$ for $i=0, 1$, $x^t(z_0A_0+z_1A_1)=0$, $(z_0A_0+z_1A_1)y=0$.
We may assume $z_1\neq 0$ and $y^tA_0y=0$. Pick $\lambda$ such that $x^tA_0x+\lambda^2y^tA_0y=0$.
By the assumptions we get $z_0(x^tA_0x+\lambda^2y^tA_0y)+z_1(x^tA_1x+\lambda^2y^tA_1y)=0$
which implies $x^tA_1x+\lambda^2y^tA_1y=0$. We get
$(x+\lambda y)^tA_i(x+\lambda y)=0$ for $i=0, 1$, as we wanted.\qed

\begin{thm0}\label{3cond}
Let $A$ be a $2\times n\times n$ matrix
with two symmetric $n\times n$ slices $A_0$, $A_1$.
The following are equivalent

(i) $Det(A)\neq 0$.

(ii) the characteristic polynomial $\det(A_0+tA_1)$ has $n$ distinct roots.

(iii) the codimension two subvariety intersection of the quadrics $A_0$, $A_1$ is smooth.
\end{thm0}

\proof  
(i) $\Longrightarrow$ (ii) Assume that $t=0$ is a double root of $\det(A_0+tA_1)=0$
We may assume that $A_0$, $A_1$ have the same shape than in the proof of Prop. \ref{4bb},
and the same argument shows that $A$ is degenerate.

(ii) $\Longrightarrow$ (iii) Assume that  the point $(1,0,\ldots,0)$ belongs to both quadrics
and their tangent spaces in this point do not meet transversely, so we may assume they are both equal to $x_1=0$.
Hence both matrices have the form

$$ A_i=\left[\begin{array}{ccc}
0&a_i&0\\
a_i&*&*\\
0&*&0\\
\end{array}\right]$$
and $\det(A_0+tA_1)$ contains the factor $(a_0+ta_1)^2$, against the assumption.

(iii) $\Longrightarrow$ (i) If $A$ is degenerate, by the Prop. \ref{degensim} we may assume that
the two slices have the same shape as in the proof of Theorem \ref{schlafli}. Hence
$A_0$ is singular at the point $(1,0,\ldots, 0)$ which is common to $A_1$,
which gives a contradiction.\qed 

\begin{prop}\label{simultaneousdiag}
Let $A$ be a $2\times n\times n$ matrix
with two symmetric $n\times n$ slices $A_0$, $A_1$.
If $Det(A)\neq 0$ then the two quadrics $A_i$ are simultaneously diagonalizable (as quadratic forms), that is there is an invertible matrix $C$
such that $C^{t}A_iC=D_i$ with $D_i$ diagonal. The columns of $C$ correspond to the $n$ distinct singular points
found for each root of $\det(A_0+tA_1)$. 
\end{prop}
\proof We may assume that $A_0$, $A_1$ are both nonsingular and from Theorem \ref{3cond}
we get distinct $\lambda_i$ for $i=1,\ldots ,n$ such that $\det(A_0+\lambda_iA_1)=0$.
For any $i=1,\ldots n$ we obtain nonzero $v_i\in\C^n$ such that  $(A_0+\lambda_iA_1)v_i=0$. From these equations
we get $\lambda_i(v_j^tA_1v_i)=v_j^tA_0v_i=\lambda_j(v_j^tA_1v_i)$.
Hence for $i\neq j$ we get $v_j^tA_1v_i=0$ and also $v_j^tA_0v_i=0$. Let $C$ be the matrix having $v_i$ as columns.
The identities found are equivalent to $C^{t}A_iC=D_i$ where $D_i$ has $v_1^tA_iv_1,\ldots,v_n^tA_iv_n$
on the diagonal . It remains to show that $v_i$ are independent.
This follows because $v_i^tA_0v_i\neq 0$ for any $i$, otherwise $A_0v_i=0$ and $A_0$ should be singular.\qed 
\begin{rem}
One may assume that
$D_0=diag(\lambda_1,\ldots,\lambda_n)$ and $D_1=diag(\mu_1,\ldots,\mu_n)$ and in this case $ Det(A)$ is proportional to
$$\prod_{i<j}(\lambda_i\mu_j-\lambda_j\mu_i)^2$$ So there are simultaneously diagonalizable pairs of quadrics
with vanishing hyperdeterminant, in other words, the converse to 
Proposition \ref{simultaneousdiag} does not hold. The condition for a pair of quadrics to be simultaneously diagonalizable is more subtle.
For two smooth conics in the plane ($n=3$) one has to avoid just the case
that the two conics touch in a single point (they can touch in two distinct points and still being simultaneously diagonalizable). 
\end{rem}

Oeding considers in \cite{Oed} the case of homogeneous polynomials of degree $d$ in $n+1$ variables,
they give a symmetric tensor of format $(n+1)\times\ldots\times (n+1)$ ($d$ times),
corresponding to the embedding $S^d\C^{n+1}\subset\otimes^{d}\C^{n+1}$.
The coefficients $a_{i_1,\ldots, i_d}$ of the multidimensional matrix satisfy
$$a_{i_1,\ldots, i_d}=a_{\sigma(i_1),\ldots, \sigma(i_d)}$$ for every permutation $\sigma$.
Don't confuse this notion with the determinantal representation like in Example \ref{333}.
For example the Fermat cubic $x_0^3+x_1^3+x_2^3$ defines a $3\times 3\times 3$ tensors
with only three entries equal to $1$, and its three determinantal representations
are all equal to $x_0x_1x_2$, corresponding to three lines.

The first easy result is the following
\begin{thm0}\label{firsteasy}
The discriminant of $f\in S^d\C^{n+1}$ divides the hyperdeterminant of 
the multidimensional matrix in $\otimes^{d}\C^{n+1}$ corresponding to $f$.
\end{thm0} 
\proof Let $f$ be singular at
$v_0$, we get that $f$ correspond to a multilinear map 
$A_f\colon\C^{n+1}\times\ldots\times\C^{n+1}\to\C$ such that
$A_f(v_0,\ldots, v_0)=0$ and
$A_f(\C^{n+1},v_0,\ldots, v_0)=A_f(v_0,\C^{n+1},v_0,\ldots, v_0)=\ldots =0$.
Hence the kernel of $A_f$ contains $v_0\otimes\ldots\otimes v_0$, $A_f$ is degenerate
and it has zero hyperdeterminant.\qed
\vskip 0.5cm

Oeding proves that the converse is true only in two cases:
the square case $n\times n$ and the $2\times 2\times 2$ case.

In all the other cases the hyperdeterminant of a symmetric tensor has the discriminant as a factor
but contains interesting extra terms.

For example, in the $3\times 3\times 3$ case the hyperdeterminant has degree 36
and it is the product
$D\cdot S^6$, where $D$ is the discriminant of degree $12$ and $S$ is the Aronhold invariant, which vanishes on plane cubics
which are sum of three cubes of linear forms. In other words $S$
is the equation of (the closure of) the $SL(3)$-orbit of the Fermat cubic $x_0^3+x_1^3+x_2^3$.

In order to describe the extra terms let's consider  any partition $\lambda=(\lambda_1,\ldots, \lambda_s)$ of $d$,
that is $d=\lambda_1+\ldots+\lambda_s$, we may assume $\lambda_1\ge\ldots\ge\lambda_s$.
For any partition $\lambda$ of $d$ define $Chow_{\lambda}(\P n)$
as the closure in $\PP (S^d\C^{n+1})$ of the set of polynomials of degree $d$ which are expressible as
$l_1^{\lambda_1}\cdots l_s^{\lambda_s}$ where $l_i$ are linear forms.
\begin{thm0}\label{oeding1}\cite{Oed} The dual variety $Chow_{\lambda}(\P n)^{\vee}$ is a hypersurface except for the two cases

(i) $n=1$ and $\lambda_s=1$

(ii) $n\ge 2$ and $\lambda=(d-1,1)$.

\end{thm0}

Let $\Theta_{\lambda,n}$ be the equation of  $Chow_{\lambda}(\P n)^{\vee}$
when it is a hypersurface (see Theorem \ref{oeding1}).

\begin{thm0}\label{oeding2}\cite{Oed}
The hyperdeterminant  of a symmetric matrix of format $n\times\ldots\times n$ ($d$ times)
splits as the product

$$\prod_{\lambda}\Theta_{\lambda,n}^{m_{\lambda}}$$
where $m_{\lambda}$ is the multinomial coefficient
${d\choose{\lambda_1,\ldots\lambda_s}}$
and the product is extended over all the partitions such that $Chow_{\lambda}(\P n)^{\vee}$ 
is a hypersurface, classified by Theorem \ref{oeding1}.
\end{thm0}

There is always the factor corresponding to the trivial partition $d$.
The Chow variety $Chow_{d}(\P n)$ is the Veronese variety and its dual variety is the discriminant
(according to Theorem \ref{firsteasy}), appearing in the product with exponent one.
Indeed one sees immediately from Theorem \ref{oeding1} that this is the only factor just in the 
square case $n\times n$ ($d=2$) and in the $2\times 2\times 2$ case.

\section{Weierstrass canonical form and Kac's Theorem}

Note that the only format $2\times b\times c$ where the hyperdeterminant exists
(so that the triangular inequality is satisfied)
are $2\times k\times k$ and $2\times k\times (k+1)$.

The $2\times k\times k$ case has the same behaviour as the symmetric case considered in Section \ref{symmetric}.
We record the main classification result in the nondegenerate case.

\begin{thm0} [Weierstrass]
Let $A$ be a $2\times k\times k$ matrix and let $A_0, A_1$ be the two slices.
Assume that $Det(A)\neq 0$. Under the action of $GL(k)\times GL(k)$ $A$ is equivalent to a matrix where
$A_0$ is the identity and $A_1=diag(\lambda_1,\ldots ,\lambda_k)$.
In this form the hyperdeterminant of $A$ is equal to $\prod_{i<j}(\lambda_i-\lambda_j)^2$.
\end{thm0}

The other case $2\times k\times (k+1)$ has boundary format and it was also solved by Weierstrass.

\begin{thm0} [Weierstrass]
\label{2kk+1}
All nondegenerate matrices of type $2\times k\times (k+1)$ are
$GL(k)\times GL(k+1)$ equivalent to the identity matrix having the two slices
$$\left[\begin{array}{cccc}1\\
&\ddots\\
&&1&\quad\\
\end{array}\right]\qquad
\left[\begin{array}{cccc}\quad&1\\
&&\ddots\\
&&&1\\
\end{array}\right]
$$
\end{thm0}

\begin{proof}
Let $A$, $A'$ two such matrices.
Since they are nondegenerate they define two exact sequences on $\PP^1$

\[0\to\O (-k)\rig{}\O^{k+1}\rig{A} \O(1)^{k}\to 0\]
\[0\to\O (-k)\rig{}\O^{k+1}\rig{A'} \O(1)^{k}\to 0\]
We want to show that there is a commutative diagram
\[\begin{array}{ccccccccc}
0&\to &\O (-k)&\rig{}&\O^{k+1}&\rig{A} &\O(1)^{k}&\to &0\\
&&\dow{1}&\searrow&\dow{f}\\
0&\to &\O (-k)&\rig{}&\O^{k+1}&\rig{A'} &\O(1)^{k}&\to &0\\
\end{array}\]
In order to show the existence of $f$ we
apply the functor $Hom(-,\O^{k+1})$ to the first row.
We get

{\footnotesize
$$Hom(\O^{k+1},\O^{k+1})\rig{g} Hom(\O(-k),\O^{k+1})\to
Ext^1(\O(1)^{k},\O^{k+1})\simeq H^1(\O(-1)^{k(k+1)})=0  $$} 
Hence $g$ is surjective and $f$ exists.
Now it is straightforward to complete the diagram with a 
 morphism $\phi\colon \O(1)^{k}\to \O(1)^{k}$, which is a 
isomorphism by the snake lemma.
\end{proof}

Let $(x_0,x_1)$ be homogeneous coordinates on $\P 1$. The identity matrix 
appearing in Theorem \ref{2kk+1} corresponds to the morphism of vector bundles
given by

\[ I_k(x_0, x_1):= \left ( \begin{array}{llll}
 x_0 & x_1 &       &        \\
     & \ddots& \ddots   &       \\
     &     & x_0&x_1 \\
     \end{array} \right ) \]

\vspace{3mm}

It is interesting, and quite unexpected, that the format
$2\times k\times (k+1)$ is a building block for all the other formats
$2\times b\times c$. The canonical form illustrated by the following
Theorem is called the Weierstrass canonical form (there is an extension 
in the degenerate case that we do not pursue here).

\begin{thm0} [Kronecker, 1890]
Let $2\le b<c$. There exist unique $n, m, q\in\N$
satisfying

$$\left\{\begin{array}{ccccc}b&=&nq&+&m(q+1)\\
c&=&n(q+1)&+&m(q+2)\\
\end{array}\right.$$ 

such that the general tensor $t\in\C^2\otimes\C^b\otimes\C^c$ decomposes 
under the action of $GL(b)\times GL(c)$ as
$n$ blocks $2\times q\times(q+1)$ and $m$ blocks $2\times (q+1)\times(q+2)$ in Weierstrass form. 
\end{thm0}

Kac has generalized this statement to the format $2\le w\le s\le t$ satisfying the inequality $t^2-wst+s^2\ge 1$.
Note that in these cases the hyperdeterminant does not exist (for $w\ge 3$).
The result is interesting because it gives again a canonical form.

Given $w$, define by the recurrence relation  
$a_0=0$, $a_1=1$, $a_j=wa_{j-1}-a_{j-2}$

For $w=2$ get $0, 1, 2,\ldots $ and Kronecker's  result.

For $w=3$ get $0, 1,3 , 8, 21, 55, \ldots $
(odd Fibonacci numbers)
\begin{figure}[h]
   \centering
\includegraphics[width=50mm]{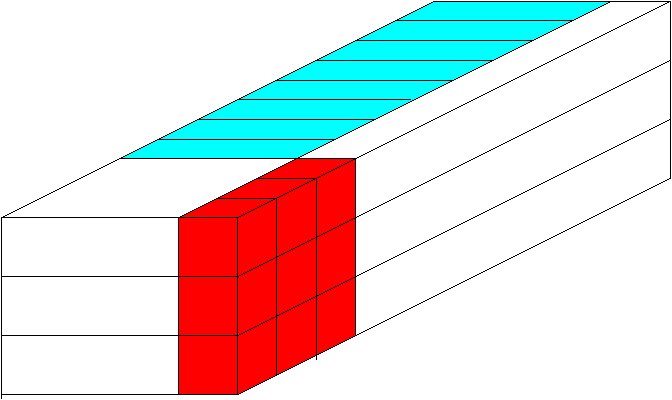}
    \caption{A decomposition in two Fibonacci blocks}
    \end{figure}

\begin{thm0} [Kac, 1980]
Let $2\le w\le s\le t$ satisfying the inequality $t^2-wst+s^2\ge 1$.
Then there exist unique $n, m, j\in\N$
satisfying

$$\left\{\begin{array}{ccccc}s&=&na_j&+&ma_{j+1}\\
t&=&na_{j+1}&+&ma_{j+2}\\
\end{array}\right.$$ 
such that the general tensor $t\in\C^w\otimes\C^s\otimes\C^t$ decomposes 
under the action of $GL(s)\times GL(t)$as
$n$ blocks $w\times a_j\times a_{j+1}$ and $m$ blocks $w\times a_{j+1}\times a_{j+2}$ which are denoted ``Fibonacci blocks''.
They can be described by representation theory (see \cite{Br}). 
\end{thm0}

The original proof of Kac (\cite{Kac}) uses representations of quivers. 
In \cite{Br} there is an independent proof in the language of vector bundles.

\section{The rank and the k-secant varieties}

The classical determinant is the equation of the dual variety of the variety of decomposable tensors.
This property of the determinant has been chosen as definition of hyperdeterminant in the multidimensional setting.

The determinant gives the condition for a homogeneous linear system to have
 nontrivial solutions.
We have discussed this second property in the multidimensional setting in Section \ref {multilinear}.

A third property of the determinant is that it vanishes precisely on matrices
not of maximal rank.

This property generalizes in the multidimensional setting in a completely different manner,
and it is no more governed by the hyperdeterminant.

The rank $r$ of a multidimensional matrix $A$ of format $(k_0+1)\times\ldots\times (k_p+1)$ 
is the minimum number of decomposable summands $t_i=x_0^i\otimes\ldots\otimes x_{n}^i$
needed to express it, that is $A=\sum_{i=1}^rt_i$ is minimal.

So the $r$-th secant variety $\sigma_r(\P {{k_0}}\times\ldots\times \P {{k_p}})$ parametrizes multidimensional
matrices of rank $r$ and their limits.

The problem of describing these secant varieties is widely open (for the first properties see \cite{CGG}). 
Even their dimension is not known in general,
although it is conjectured that it coincides with the expected value $r(1+\sum_{i=1}^rk_i)-1$
(when this number is smaller than the dimension of the ambient space) unless a list of well understood cases
\cite{AOP}. The analogous result for symmetric multidimensional matrices has been proved by Alexander and Hirschowitz.

The first attempt to find equations of $\sigma_r(\P {{k_0}}\times\ldots\times \P {{k_p}})$
is through the minors of the flattening maps defined in Proposition \ref{rankcontraction}.

Indeed Raicu proved in \cite{Ra} that the ideal of $\sigma_2(\P {{k_0}}\times\ldots\times \P {{k_p}})$
is generated by the $3\times 3$-minors of $C_i(\phi)$, so proving a conjecture by Garcia, Stillmann and Sturmfels. 

These varieties $\sigma_2$ are never hypersurfaces, unless the trivial case of $\sigma_2(\P {{2}}\times \P {{2}})$,
which is given by the classical $3\times 3$ determinant.

The first nontrivial case is given by 
 $\sigma_i(\P {{2}}\times \P {{2}}\times \P {{2}})$,
which can be described in a uniform way for $i=1,\ldots 4$ by the $2i+1$-minors
of the more general flattening (Young flattening)

$$C_0(\phi)\colon V_0^{\vee}\otimes V_1\to \wedge^2V_1\otimes V_2$$

for $\phi\in V_0\otimes V_1\otimes V_2$
and by the $2i+1$-minors of the analogous flattening 
$C_1(\phi)$ and $C_2(\phi)$ obtained by permutations.

For $i=4$ the secant variety $\sigma_4(\P {{2}}\times \P {{2}}\times \P {{2}})$ is a hypersurface, in this case
$\det C_i(\phi)$ are independent from $i$ and they all give the same hypersurface of degree $9$
(see \cite{Ot}\S 3). This hypersurface
was found first by Strassen in \cite{St}, with a slightly different construction.
We emphasize that this hypersurface is different from the hyperdeterminant, which has degree $36$, constructed in  
Example \ref{333}.

If $A_i$ are the three slices of $A$, the matrix of $C_0(\phi)$ can be depicted as

$$ \left[\begin{array}{rrr}
0&A_2&-A_1\\
-A_2&0&A_0\\
A_1&-A_0&0\\
\end{array}\right]$$

When $A$ is symmetric, $C_i(\phi)$ are independent from $i$ and appear to be skew-symmetric.
In this case the pfaffians of order $2i+2$ of $C_0(\phi)$ define the $i$-th secant variety
of the Veronese variety given by the $3$-embedding of $\P 2$ in $\PP (S^3\C^3)$.
For $i=3$ we get the Aronhold invariant of degree $4$ which is the equation of the orbit of the
Fermat cubic (see \cite{LO}), that we encountered in Section \ref{symmetric}.

\section{Open problems }

There are a lot of interesting (and difficult) open problems 
on the subject, starting from  looking for equations of secant varieties to the Segre varieties.

Here we propose three problems on the hyperdeterminant that seem tractable (at least the first two) and interesting to me .
\vskip 0.5cm

{\bf Problem 1} Find the equations for the dual varieties to Segre varieties when they are not hypersurfaces,
so when $k_0=\max_j k_j>\sum k_i$ (when they are hypersurfaces the single equation is the hyperdeterminant).

Let's see the example of format $4\times 2\times 2$.
In this case the dual variety to $\P 3\times\P 1\times\P 1$ has codimension $2$,
by Theorem \ref{dimdual} (i).
By Theorem \ref{condi0} the dual variety consists of matrices which are not $0$-degenerate.
One sees that if $A$ has format $4\times 2\times 2$,
the multilinear system $\tilde A\colon\C^2\otimes\C^2\to\C^4$ has a nontrivial solution if and only if 
the following two conditions hold

i) the hyperdeterminant of every 
$3\times 2\times 2$ submatrix of $A$ vanishes. 

ii) $\det(\tilde A)=0$ where $\tilde A$ is seen as a $4\times 4$ matrix.

So the equations in i) and ii) give the answer to this problem for the format $4\times 2\times 2$.

Note that the equations of the dual of $\P{{k_0}}\times\ldots\times\P{{k_p}}$ when
$k_0>\sum k_i$ must contain all the hyperdeterminants
of submatrices of boundary format (with $k_0'=\sum k_i$) .
\vskip 0.5cm

{\bf Problem 2} Compute the irreducible
factors of the hyperdeterminant of a skew-symmetric tensor  in $\wedge^{d}\C^n\subset\otimes^d\C^n$.

This means to extend Oeding Theorem \ref{oeding2} to the skew-symmetric case.
In this case even in the square case the classical determinant is not irreducible, indeed it is the square of the pfaffian.
The dual varieties to Grassmann varieties will play into the game.
\vskip 0.5cm

{\bf Problem 3}
  This question is a bit more vague. The definition of hyperdeterminant of boundary format with
the linear map $\partial_A$ (compare with Theorem \ref{main}
and Corollary \ref{defdet}) can be generalized to
  other cases where the codimension of the degenerate matrices is
bigger than one. 
Specifically, if $k_0, \ldots, k_p$ are
  nonnegative integers satisfying $k_0=\sum_{i=1}^pk_i$ then we denote
  again $m_j=\sum_{i=1}^{j-1}k_i$ with the convention $m_1=0$, like in Section \ref{multilinear}.  
  
  Assume we have vector spaces $V_0,\ldots ,V_p$ and a positive
  integer $q$ such that $\dim V_0=q(k_0+1), \dim V_1=q(k_1+1)$ and $\dim
  V_i=(k_i+1)$ for $i=2,\ldots ,p$.  Then the vector spaces
  $V_0\otimes S^{m_1}V_1\otimes\ldots\otimes S^{m_p}V_p$ and $
  S^{m_1+1}V_1\otimes\ldots\otimes S^{m_p+1}V_p$ still have the same
  dimension and there is an analogous invariant given by $\det(\partial_A)$ .
The question is to study the properties of this invariant. 
  
  Although this construction can seem artificious, it found an application in the first case  $q=p=2$,
  leading to the proof \cite{CO}
  that the moduli space of instanton bundles on $\PP^3$
is affine.

\end{document}